\documentclass[11pt,fleqn]{article}
\usepackage[utf8]{inputenc}
\usepackage[top=1in,bottom=1in,left=1in,right=1in]{geometry}
\usepackage{fouriernc}
\usepackage{authblk}
\usepackage{mathtools}
\usepackage[all,cmtip]{xy}
\usepackage{multicol}
\usepackage{amssymb}
\usepackage{amsthm}
\usepackage[all]{xy}

\newtheorem{theorem}{Theorem}

\theoremstyle{definition}

\theoremstyle{remark}
\newtheorem{remark}[theorem]{Remark}
\title{The equitable presentation of $\mathfrak{osp}_q(1|2)$\\ and a $q$-analog of the Bannai--Ito algebra}
\author[1]{Vincent X. Genest}
\author[1]{Luc Vinet}
\author[2]{Alexei Zhedanov}
\affil[1]{Centre de recherches math\'ematiques, Universit\'e de Montr\'eal, P.O. Box 6128, Centre-ville Station, Montr\'eal, Canada H3C 3J7}
\affil[2]{Donetsk Institute for Physics and Technology, Donetsk 83114, Ukraine}

\date{}
\numberwithin{equation}{section}

\begin{document}
\maketitle
\thispagestyle{empty}
\hrule
\begin{abstract}
\noindent
The equitable presentation of the quantum superalgebra $\mathfrak{osp}_q(1|2)$, in which all generators appear on an equal footing, is exhibited. It is observed that in their equitable presentations, the quantum algebras $\mathfrak{osp}_q(1|2)$ and $\mathfrak{sl}_q(2)$ are related to one another by the formal transformation $q\rightarrow -q$. A $q$-analog of the Bannai--Ito algebra is shown to arise as the covariance algebra of $\mathfrak{osp}_q(1|2)$.
\smallskip

\noindent \textbf{AMS classification numbers:} 17B37, 20G42, 81R50
\end{abstract}
\hrule
\section{Introduction}
The purpose of this Letter is threefold: to display the equitable, or $\mathbb{Z}_3$-symmetric, presentation of the quantum superalgebra $\mathfrak{osp}_q(1|2)$, to show that the equitable presentations of $\mathfrak{osp}_q(1|2)$ and $sl_{q}(2)$ are related to one another by the formal transformation $q\rightarrow -q$, and to demonstrate that the covariance algebra of $\mathfrak{osp}_q(1|2)$ is a $q$-analog of the Bannai--Ito algebra.

Our considerations take root in the Racah problem for the $\mathfrak{su}(2)$ algebra, i.e. the coupling of three angular momenta. In this problem, the states are usually described in terms of the quantum numbers $j_i$ associated to the individual angular momenta $\vec{J}_i$ with $i\in \{1,2,3\}$, the quantum number $j$ associated to the total angular momentum $\vec{J}=\vec{J}_1+\vec{J}_2+\vec{J}_3$, the quantum number $M$ associated to the projection of the total angular momentum $\vec{J}$ along one axis, and any one of the quantum numbers $j_{12}$, $j_{23}$, $j_{31}$ associated to the intermediate angular momenta $\vec{J}_{ij}=\vec{J}_i+\vec{J}_j$ for $(ij)\in \{(12), (23), (31)\}$. These bases are related via Racah coefficients \cite{1996_Edmonds}. The main drawback of such bases is their involved behavior under particle permutations. To circumvent this problem, Chakrabarti \cite{1964_Chakrabarti_AnnIHP_1_301}, L\'evy-Leblond and L\'evy-Nahas \cite{1965_LevyLeblond&LevyNahas_JMathPhys_6_1372} devised an ``equitable'' coupling scheme and showed that there is a ``democratic'' basis specified by the quantum numbers $j_1$, $j_2$, $j_3$, $j$ and $\zeta$, where $\zeta$ is the eigenvalue of the volume operator $\Delta=(\vec{J}_1\times \vec{J}_2)\cdot \vec{J}_3$. The three angular momenta $\vec{J}_i$ enter symmetrically in this scheme and the states of the democratic basis have definite behaviors under particle permutations.

The Racah--Wilson algebra is the hidden algebraic structure behind the Racah problems of $\mathfrak{su}(2)$ and $\mathfrak{su}(1,1)$ \cite{1988_Granovskii&Zhedanov_JETP_94_49}. The concern for a democratic approach to these Racah problems leads to the equitable presentation of this algebra \cite{2014_Genest&Vinet&Zhedanov_JPhysA_47_025203}. In this presentation, the defining relations of the Racah--Wilson algebra are $\mathbb{Z}_3$-symmetric and all the generators appear on an equal footing, whence the epithets ``equitable'' or ``democratic''. It was recently shown in \cite{2013_Gao&Wang&Hou_LinAlgAppl_439_1834} that the equitable generators of the Racah/Wilson algebra can also be realized as quadratic expressions in the equitable $\mathfrak{sl}(2)$ generators proposed in \cite{2007_Hartwig&Terwilliger_JAlg_840}. Note that the Racah--Wilson algebra also arises as symmetry algebra for superintegrable systems \cite{2007_Kalnins&Miller&Post_JPhysA_40_11525} and encodes the bispectrality of the Racah/Wilson polynomials \cite{2014_Genest&Vinet&Zhedanov_LettMathPhys_104_931}.

The Racah problem can also be posited for the quantum algebra $\mathfrak{sl}_q(2)$. In this case, it is the Askey--Wilson algebra \cite{1991_Zhedanov_TheorMathPhys_89_1146}, also known as the Zhedanov algebra \cite{2007_Koornwinder_SIGMA_3_63}, that arises as the hidden algebraic structure \cite{1993_Granovskii&Zhedanov_JGroupThPhys_1_161}. This algebra encodes the bispectrality of the Askey-Wilson polynomials and, as shown in \cite{1993_Granovskii&Zhedanov_JPhysA_26_L357}, arises as the covariance algebra for $\mathfrak{sl}_q(2)$. An equitable presentation of the (universal) Askey-Wilson algebra was offered by Terwilliger in \cite{2011_Terwilliger_SIGMA_7_99} who also showed that it can be realized by quadratic combinations of the equitable generators of $\mathfrak{sl}_q(2)$. The equitable presentation of $\mathfrak{sl}_q(2)$ was itself studied in \cite{2006_Ito&Terwilliger&Weng_JAlg_298_284}. A democratic presentation for the quantum group $U_{q}(\mathfrak{g})$ associated with a symmetrizable Kac-Moody algebra $\mathfrak{g}$ was also  proposed in \cite{2006_Terwilliger_JAlg_298_302}.

In a recent paper \cite{2015_Genest&Vinet&Zhedanov_ArXiv_150105602}, the Racah problem for the quantum superalgebra $\mathfrak{osp}_q(1|2)$ was considered. It was shown that in this case a $q$-analog of the Bannai--Ito algebra, an algebra proposed in \cite{2012_Tsujimoto&Vinet&Zhedanov_AdvMath_229_2123}, appears as the ``hidden'' algebraic structure. The algebra obtained in \cite{2015_Genest&Vinet&Zhedanov_ArXiv_150105602} exhibits a $\mathbb{Z}_3$ symmetry and is related to the Askey--Wilson algebra by the formal transformation $q\rightarrow -q$.

In this Letter, we display the equitable presentation of the quantum superalgebra $\mathfrak{osp}_q(1|2)$, determine its relation with the equitable presentation of $\mathfrak{sl}_q(2)$ and show that it can be used to realize the $q$-analog of the Bannai--Ito algebra defined in \cite{2015_Genest&Vinet&Zhedanov_ArXiv_150105602}. The results that we present here enrich the understanding of the quintessential quantum superalgebra $\mathfrak{osp}_q(1|2)$ and shed light on its relationship with other algebraic structures that have appeared recently. The contents of the Letter are as follows.

In Section 2, the definition of $\mathfrak{osp}_q(1|2)$ is reviewed and its extension by the grade involution is presented. A two-parameter family of $\mathfrak{osp}_q(1|2)$-modules is defined. The equitable presentation of $\mathfrak{osp}_q(1|2)$ is introduced and several expressions are given for the Casimir operator. The equitable presentation of $\mathfrak{sl}_q(2)$ is reviewed and compared with the one found for $\mathfrak{osp}_q(1|2)$. In section 3, the realization of the $q$-deformed Bannai--Ito algebra in terms of the equitable $\mathfrak{osp}_q(1|2)$ generators is presented. A short conclusion follows.
\section{The $\mathfrak{osp}_q(1|2)$ algebra and its equitable presentation}
In this section, we recall the definition of the $\mathfrak{osp}_q(1|2)$ algebra, present its extension by the grade involution, define a family of irreducible representations and display its equitable presentation. We review the equitable presentation of $sl_{q}(2)$ and compare it with that of $\mathfrak{osp}_q(1|2)$.
\subsection{Definition of $\mathfrak{osp}_q(1|2)$, the grade involution, and representations}
Let $q$ be a complex number which is not a root of unity and let $[n]_{q}$ denote
\begin{align*}
 [n]_{q}=\frac{q^{n}-q^{-n}}{q-q^{-1}}.
\end{align*}
The quantum superalgebra $\mathfrak{osp}_q(1|2)$ is the $\mathbb{Z}_2$-graded unital associative $\mathbb{C}$-algebra generated by the even element $A_0$ and the odd elements $A_{\pm}$ satisfying the relations \cite{1989_Kulish&Reshetikhin_LettMathPhys_18_143}
\begin{align*}
 [A_0,A_{\pm}]=\pm A_{\pm},\qquad \{A_{+},A_{-}\}=[2A_0]_{q^{1/2}},
\end{align*}
where $[x,y]=xy-yx$ and $\{x,y\}=xy+yx$ respectively stand for the commutator and the anticommutator. The sCasimir operator of $\mathfrak{osp}_q(1|2)$ is defined as \cite{1995_Lesniewski_JMathPhys_36_1457}
\begin{align*}
 S=A_{+}A_{-}-[A_0-1/2]_{q}.
\end{align*}
This operator is readily seen to obey the relations
\begin{align*}
 \{S, A_{\pm}\}=0,\qquad [S,A_0]=0.
\end{align*}
The abstract $\mathbb{Z}_2$-grading of $\mathfrak{osp}_q(1|2)$ can be concretized by adding the grade involution $P$ to the set of generators and by declaring that the even and odd generators respectively commute and anticommute with $P$. The quantum superalgebra $\mathfrak{osp}_q(1|2)$ can thus be presented as the unital associative  $\mathbb{C}$-algebra generated by the elements $A_0$, $A_{\pm}$ and the involution $P$ obeying the relations
\begin{subequations}
\label{OSP}
\begin{align}
\label{OSPA}
[A_0,A_{\pm}]=\pm A_{\pm},\quad \{A_{+},A_{-}\}=[2A_0]_{q^{1/2}},\quad
 \{P,A_{\pm}\}=0,\quad [P,A_0]=0,\quad  P^2=1.
\end{align}
In \eqref{OSPA}, the parity of the elements no longer needs to be specified. Upon introducing the generators 
\begin{align*}
 K=q^{A_0},\qquad K^{-1}=q^{-A_0},
\end{align*}
one can write the relations \eqref{OSPA} in the form
\begin{align}
 \begin{aligned}
  &K A_{+} K^{-1}=q A_{+},\qquad K A_{-}K^{-1}=q^{-1}A_{-},\qquad KK^{-1}=1,\quad P^2=1,
  \\
  &[K,P]=0,\qquad [K^{-1},P]=0,\qquad \{A_{\pm}, P\}=0,\qquad \{A_{+},A_{-}\}=\frac{K-K^{-1}}{q^{1/2}-q^{-1/2}}.
 \end{aligned}
\end{align}
\end{subequations}
It is directly verified that the Casimir operator
\begin{align}
\label{Cas}
 Q=\Big(A_{+}A_{-}-[A_0-1/2]_{q}\Big)P,
\end{align}
which is related to the sCasimir operator by $Q=SP$, commutes with all the generators in \eqref{OSP}.  The quantum superalgebra $\mathfrak{osp}_q(1|2)$ can be equipped with a Hopf algebraic structure \cite{1989_Kulish&Reshetikhin_LettMathPhys_18_143}. Introduce the coproduct map $\Delta:\mathfrak{osp}_q(1|2)\rightarrow \mathfrak{osp}_{q}(1|2)\otimes \mathfrak{osp}_q(1|2)$ with
\begin{align}
\label{Delta}
\begin{aligned}
 &\Delta(A_{+})=A_{+}\otimes KP+1\otimes A_{+},\qquad \Delta(A_{-})=A_{-}\otimes P+K^{-1}\otimes A_{-},
 \\
 &\Delta(K)=K\otimes K,\quad \Delta(P)=P\otimes P,
\end{aligned}
 \end{align}
the counit map $\epsilon:\mathfrak{osp}_q(1|2)\rightarrow \mathbb{C}$ with
\begin{align}
\label{Epsilon}
 \epsilon(P)=1,\quad \epsilon(K)=1,\quad \epsilon(A_{\pm})=0,
\end{align}
and the coinverse map $\sigma:\mathfrak{osp}_q(1|2)\rightarrow \mathfrak{osp}_q(1|2)$  with
\begin{align}
\label{Sigma}
 \sigma(P)=P,\quad \sigma(K)=K^{-1},\quad \sigma(A_{+})=-A_{+}K^{-1}P,\quad \sigma(A_{-})=-K A_{-}P.
\end{align}
It is a straightforward exercise to verify that with the coproduct $\Delta$, the counit $\epsilon$ and the coinverse $\sigma$ defined as in \eqref{Delta}, \eqref{Epsilon} and \eqref{Sigma}, the algebra \eqref{OSP} complies with the well-known requirements for a Hopf algebra \cite{2011_Underwood}. 

\begin{remark}
Let us note that in computing with the coproduct, one should use the standard product rule $(a\otimes b)(c\otimes d)=(ac\otimes bd)$, as opposed to the usual graded product rule used for superalgebras when the grade involution is not introduced.
\end{remark}

Let us now bring a two-parameter family of irreducible representations of $\mathfrak{osp}_q(1|2)$. Let $\nu$ be a real parameter and let $e=\pm 1$. Moreover, let $W^{(e,\nu)}$ be the infinite-dimensional  vector space spanned by the basis vectors $f_{n}^{(e,\nu)}$, where $n$ is a non-negative integer. It is verified that the actions
\begin{alignat}{2}
\label{Rep}
\begin{aligned}
 K\,f_{n}^{(e,\nu)}&=q^{n+\nu+1/2}\,f_{n}^{(e,\nu)},\quad &P\, f_{n}^{(e,\nu)}&=e (-1)^{n}\,f_{n}^{(e,\nu)},
 \\
 A_{+}\,f_{n}^{(e,\nu)}&= \,f_{n+1}^{(e,\nu)},\quad& A_{-}\,f_{n}^{(e,\nu)}&=\rho_{n} \,f_{n-1}^{(e,\nu)},
 \end{aligned}
\end{alignat}
where
\begin{align*}
 \rho_n=[n+\nu]_{q}-(-1)^{n}[\nu]_{q},
\end{align*}
define representations of $\mathfrak{osp}_q(1|2)$ on $W^{(e,\nu)}$. For generic values of $\nu$, one has $\rho_n \neq 0$ for all $n\geqslant 1$. As a consequence, these representations are irreducible. On $W^{(e,\nu)}$, the Casimir operator \eqref{Cas} acts as a multiple of the identity
\begin{align*}
 Q\,f_{n}^{(e,\nu)}=-\,e\,[\nu]_{q}\,f_{n}^{(e,\nu)},
\end{align*}
as expected from Schur's lemma. Note that the representations $W^{(e,\nu)}$ are associated to the $q$-analog of the parabosonic oscillator \cite{1990_Floreanini&Vinet_JPhysA_23_L1019}.
\begin{remark}
 It is possible to define finite-dimensional representations of \eqref{OSP}. Indeed, if one takes $\nu=-(N+1)/2$, where $N$ is a even integer, one can use the actions \eqref{Rep} to define $(N+1)$-dimensional irreducible representations of $\mathfrak{osp}_q(1|2)$.
\end{remark}
The representations $W^{(e,\nu)}$ have a Bargmann realization on functions of the complex argument $z$. In this realization, the basis vectors $f_{n}^{(e,\nu)}\equiv f_{n}^{(e,\nu)}(z)$ read
\begin{align*}
 f_{n}^{(e,\nu)}(z)=z^{n},\qquad n=0,1,2,\ldots,
 \end{align*}
and the generators have the expressions
\begin{align}
\label{Bargmann}
\begin{aligned}
 &K(z)=q^{\nu+1/2}\,T_{q},\quad P(z)=e R_{z},
 \\
 &A_{+}(z)=z,\quad  A_{-}(z)=\frac{q^{\nu}}{q-q^{-1}}\left(\frac{T_{q}-R_{z}}{z}\right)-\frac{q^{-\nu}}{q-q^{-1}}\left(\frac{T_{q}^{-1}-R_{z}}{z}\right),
 \end{aligned}
\end{align}
where $R_z g(z)=g(-z)$ is the reflection operator and where $T_{q}^{\pm} g(z)=g(q^{\pm 1}z)$ is the $q$-shift operator.
\subsection{The equitable presentation of $\mathfrak{osp}_q(1|2)$}
Let $X$, $Y^{\pm}$, $Z$ and $\omega_{y}$ be defined as
\begin{align}
\label{XYZ}
 X=K^{-1}P-(1-q^{-1}) A_{+}K^{-1}P,\quad Y^{\pm}=K^{\pm}P,\quad Z=K^{-1}P+(q^{1/2}+q^{-1/2})A_{-}P,\quad \omega_y=P.
\end{align}
Using the relations \eqref{OSP}, one readily verifies that these operators satisfy
\begin{align}
\label{Equitable-1}
 \frac{q^{1/2}XY+q^{-1/2}YX}{q^{1/2}+q^{-1/2}}=1,\qquad \frac{q^{1/2}YZ+q^{-1/2} ZY}{q^{1/2}+q^{-1/2}}=1,\qquad \frac{q^{1/2}ZX+q^{-1/2} XZ}{q^{1/2}+q^{-1/2}}=1.
\end{align}
In addition to $YY^{-1}=1$ and $\omega_y^2=1$, one has also the relations
\begin{align}
\label{Equitable-2}
 X\omega_y+\omega_yX=2\,Y^{-1}\omega_y,\qquad Y\omega_y+\omega_yY=2\,Y\omega_y,\qquad Z\omega_y+\omega_yZ=2\, Y^{-1}\omega_y.
\end{align}
We refer to the relations \eqref{Equitable-1} and \eqref{Equitable-2} as the \emph{equitable} presentation of $\mathfrak{osp}_q(1|2)$ and to the generators $X$ ,$Y^{\pm}$, $Z$ and $\omega_y$ as the equitable generators. It is observed that in this presentation, the generators are more or less on an equal footing; some asymmetry occurs in the relations with the involution $\omega_y$ given in \eqref{Equitable-2}. The standard generators of $\mathfrak{osp}_q(1|2)$ can be expressed as follows in terms of the equitable generators:
\begin{align}
 A_{+}=\frac{1-XY}{1-q^{-1}},\quad A_{-}=\frac{(Z-Y^{-1})\omega_y}{q^{1/2}+q^{-1/2}},\quad K^{\pm}=Y^{\pm}\omega_y,\quad P=\omega_y.
\end{align}
In the equitable presentation, the ``normalized'' Casimir operator
\begin{align}
\label{Norm-Cas}
 \Upsilon=(q-q^{-1})\,Q,
\end{align}
can be written in several different ways. One has
\begin{alignat*}{2}
 \Upsilon&=q^{1/2}X-q^{-1/2}Y+q^{1/2} Z-q^{1/2}XYZ,\qquad  &\Upsilon&=q^{1/2}Y-q^{-1/2}Z+q^{1/2} X-q^{1/2}YZX,
 \\
 \Upsilon&=q^{1/2}Z-q^{-1/2}X+q^{1/2}Y-q^{1/2}ZXY,\qquad& \Upsilon&=q^{1/2}Y-q^{-1/2}Z-q^{-1/2}X+q^{-1/2}ZYX,
 \\
 \Upsilon&=q^{1/2}Z-q^{-1/2}X-q^{-1/2}Y+q^{-1/2}XZY,\qquad &\Upsilon&=q^{1/2}X-q^{-1/2}Y-q^{-1/2}Z+q^{-1/2}YXZ.
\end{alignat*}
With respect to the presentation \eqref{XYZ}, the coproduct \eqref{Delta} has the expression
\begin{alignat*}{2}
 &\Delta(X)=X\otimes 1+Y^{-1}\otimes (X-1),\qquad& \Delta(Z)&=Z\otimes 1+Y^{-1}\otimes (Z-1),
 \\
 &\Delta(Y)=Y\otimes Y,\qquad &\Delta(\omega_y)&=\omega_y\otimes \omega_y.
\end{alignat*}
On the basis $f_{n}^{(e,\nu)}$, the equitable generators have the actions
\begin{align}
\label{Actions-2}
 \begin{aligned}
 X\,f_{n}^{(e,\nu)}&=e (-1)^{n}q^{-(n+\nu+1/2)}\Big(f_{n}^{(e,\nu)}-(1-q^{-1}) f_{n+1}^{(e,\nu)}\Big),
 \\
 Y\, f_{n}^{(e,\nu)}&=e(-1)^{n}q^{n+\nu+1/2}\,f_{n}^{(e,\nu)},
 \\
 Z\,f_{n}^{(e,\nu)}&=e(-1)^{n}\Big(q^{-(n+\nu+1/2)}\,f_{n}^{(e,\nu)}+(q^{1/2}+q^{-1/2})\,\rho_n \,f_{n-1}^{(e,\nu)}\Big).
 \end{aligned}
\end{align}
\subsection{The equitable presentation of $\mathfrak{sl}_q(2)$}
Let us now establish the relation between the equitable presentations of $\mathfrak{osp}_q(1|2)$ and $\mathfrak{sl}_q(2)$. The quantum algebra $\mathfrak{sl}_q(2)$ is defined as the unital $\mathbb{C}$-algebra with generators $\kappa^{\pm}$, $J_{+}$, $J_{-}$ and relations
\begin{align}
\label{slq}
 \kappa \kappa^{-1}=\kappa^{-1}\kappa=1,\quad \kappa J_{+}\kappa^{-1}= q J_{+},\quad \kappa J_{-}\kappa^{-1}=q^{-1}J_{-},\quad [J_{+},J_{-}]=\frac{\kappa-\kappa^{-1}}{q^{1/2}-q^{-1/2}}.
\end{align}
The equitable generators $x$, $y^{\pm}$ and $z$ of $\mathfrak{sl}_q(2)$ are given by \cite{2006_Ito&Terwilliger&Weng_JAlg_298_284}
\begin{align}
\label{slq-equitable}
 x=\kappa^{-1}-(q^{1/2}-q^{-1/2})J_{+}\kappa^{-1},\qquad y^{\pm}=\kappa^{\pm},\qquad z=\kappa^{-1}+(1-q^{-1})J_{-},
\end{align}
and satisfy the relations
\begin{align}
\label{Equitable-su}
 \frac{q^{1/2}xy-q^{-1/2}yx}{q^{1/2}-q^{-1/2}}=1,\qquad  \frac{q^{1/2}yz-q^{-1/2}zy}{q^{1/2}-q^{-1/2}}=1,\qquad 
  \frac{q^{1/2}zx-q^{-1/2}xz}{q^{1/2}-q^{-1/2}}=1.
\end{align}
It is directly seen that the equitable presentation of $\mathfrak{sl}_q(2)$ given in \eqref{Equitable-su} and the equitable presentation of $\mathfrak{osp}_q(1|2)$ given in \eqref{Equitable-1} are related to one another by the formal transformation $q\rightarrow -q$. This formal relation can also be observed using the standard presentations \eqref{OSP} and \eqref{slq}. Indeed, upon defining the generators
\begin{align*}
 \widetilde{\kappa}=KP,\qquad \widetilde{\kappa}^{-1}=K^{-1}P,\qquad \widetilde{J}_{+}=\frac{1}{i}\left(\frac{1-q^{-1}}{q^{1/2}+q^{-1/2}}\right)A_{+},\qquad \widetilde{J}_{-}=\left(\frac{q^{1/2}+q^{-1/2}}{1+q^{-1}}\right)A_{-}P,
\end{align*}
one finds that they satisfy the relations
\begin{align*}
 \widetilde{\kappa}\widetilde{\kappa}^{-1}=\widetilde{\kappa}^{-1}\widetilde{\kappa}=1,\quad \widetilde{\kappa} \widetilde{J}_{+}\widetilde{\kappa}^{-1}=-q \widetilde{J}_{+},\quad  \widetilde{\kappa} \widetilde{J}_{-}\widetilde{\kappa}^{-1}=-q^{-1} \widetilde{J}_{-},\quad [\widetilde{J}_{+},\widetilde{J}_-]=\frac{\widetilde{\kappa}-\widetilde{\kappa}^{-1}}{i(q^{1/2}+q^{-1/2})},
\end{align*}
which indeed corresponds to \eqref{slq} with $q\rightarrow -q$.
\begin{remark}
 Note that if one artificially introduces an involution $\omega_y$ with $\{J_{\pm},\omega_y\}=0$ and $[\kappa,\omega_y]=0$, relations of the form \eqref{Equitable-2} also appear in the equitable presentation of $\mathfrak{sl}_{q}(2)$.
\end{remark}
\section{A $q$-generalization of the Bannai--Ito algebra and the covariance algebra of $\mathfrak{osp}_q(1|2)$}
In this section, the definitions of the Bannai--Ito algebra and that of its $q$-extension are reviewed. It is shown that the $\mathbb{Z}_3$-symmetric $q$-deformed Bannai--Ito algebra can be realized in terms of the equitable $\mathfrak{osp}_q(1|2)$ generators.
\subsection{The Bannai--Ito algebra and its $q$-extension}
The Bannai--Ito algebra first arose in \cite{2012_Tsujimoto&Vinet&Zhedanov_AdvMath_229_2123} as the algebraic structure encoding the bispectral properties of the Bannai--Ito polynomials. It also appears as the hidden algebra behind the Racah problem for the Lie superalgebra $\mathfrak{osp}(1|2)$ \cite{2014_Genest&Vinet&Zhedanov_ProcAmMathSoc_142_1545} and as a symmetry algebra for superintegrable systems \cite{2015_DeBie&Genest&Vinet_ArXiv_1501.03108, 2014_Genest&Vinet&Zhedanov_JPhysA_47_205202}. The Bannai--Ito algebra is unital associative algebra over $\mathbb{C}$ with generators $K_1$, $K_2$, $K_3$ and relations
\begin{align}
\label{BI}
 \{K_1,K_2\}=K_3+\alpha_3,\quad \{K_2,K_3\}=K_1+\alpha_1,\quad \{K_3,K_1\}=K_2+\alpha_2,
\end{align}
where $\alpha_i$, $i=1,2,3$, are structure constants. This algebra admits the Casimir operator
\begin{align}
\label{BI-Cas}
 L=K_1^2+K_2^2+K_3^2,
\end{align}
which commutes with every generator $K_i$, $i=1,2,3$. In \cite{2015_Genest&Vinet&Zhedanov_ArXiv_150105602}, a $q$-deformation of the algebra \eqref{BI} was identified in the study of the Racah problem for $\mathfrak{osp}_q(1|2)$.  This $q$-extension has generators $I_1$, $I_2$, $I_3$ and relations
\begin{align}
\label{qBI}
 \{I_1, I_2\}_{q}=I_3+\iota_3,\quad \{I_2,I_3\}_{q}=I_1+\iota_1,\quad \{I_3,I_1\}_{q}=I_2+\iota_2,
\end{align}
where  $\iota_1$, $\iota_2$, $\iota_3$ are structure constants and where
\begin{align}
 \{A,B\}_{q}=q^{1/2}AB+q^{-1/2}BA,
\end{align}
is the $q$-anticommutator. The algebra \eqref{qBI} is formally related to the Zhedanov algebra by the transformation $q\rightarrow -q$ \cite{2015_Genest&Vinet&Zhedanov_ArXiv_150105602}. It has for Casimir operator
\begin{align}
\label{qBi-Cas}
 \Lambda=(q^{-1/2}-q^{3/2}) I_1I_2I_3+q I_1^2+q^{-1}I_2^2+q I_3^2-(1-q)\,\iota_1\, I_1-(1-q^{-1})\,\iota_2\,I_2-(1-q)\,\iota_3\, I_3,
\end{align}
which commutes with all generators $I_i$, $i=1,2,3$. It is easily seen that in the limit $q\rightarrow 1$, the relations \eqref{qBI} and the expression \eqref{qBi-Cas} tend to the relations \eqref{BI} and to the relation \eqref{BI-Cas}.
\subsection{Covariance algebra of $\mathfrak{osp}_q(1|2)$}
Let $a^{\pm}$, $b^{\pm}$, $c^{\pm}$ be complex parameters and consider the operators $A$, $B$, $C$ defined as follows:
\begin{align}
\label{Realization-2}
\begin{aligned}
 &A=a X-a^{-1} Y+\frac{b c^{-1}\left(XY-YX\right)}{q^{1/2}+q^{-1/2}},\qquad 
 B=b Y-b^{-1}Z+\frac{c a^{-1}\left(YZ-ZY\right)}{q^{1/2}+q^{-1/2}},
 \\
 &C=c Z-c^{-1}X+\frac{a b^{-1}\left(ZX-XZ\right)}{q^{1/2}+q^{-1/2}},
 \end{aligned}
\end{align}
where $X$, $Y$ and $Z$ are the equitable generators of $\mathfrak{osp}_q(1|2)$ defined in \eqref{XYZ}. A direct calculation shows that the elements $A$, $B$ and $C$ satisfy the relations
\begin{align}
\label{qBI-1}
\begin{aligned}
 \frac{q^{1/2} AB+q^{-1/2} BA}{q-q^{-1}}&=C+\frac{(a-a^{-1})(b-b^{-1})-(c-c^{-1})\Upsilon}{q^{1/2}-q^{-1/2}},
 \\
  \frac{q^{1/2} BC+q^{-1/2} CB}{q-q^{-1}}&=A+\frac{(b-b^{-1})(c-c^{-1})-(a-a^{-1})\Upsilon}{q^{1/2}-q^{-1/2}},
  \\
    \frac{q^{1/2} CA+q^{-1/2} AC}{q-q^{-1}}&=B+\frac{(c-c^{-1})(a-a^{-1})-(b-b^{-1})\Upsilon}{q^{1/2}-q^{-1/2}},
    \end{aligned}
\end{align}
where $\Upsilon$ is the normalized Casimir operator \eqref{Norm-Cas} of $\mathfrak{osp}_q(1|2)$. The algebra \eqref{qBI-1} is easily seen to be equivalent to the $q$-deformed Bannai--Ito algebra \eqref{qBI}. Indeed, upon taking
\begin{align}
\label{Realization-1}
 M_1=\frac{A}{q-q^{-1}},\quad M_2=\frac{B}{q-q^{-1}},\quad M_3=\frac{C}{q-q^{-1}},
\end{align}
one finds that the elements $M_1$, $M_2$, $M_3$ satisfy the defining relations \eqref{qBI} of the $q$-deformed Bannai--Ito algebra 
\begin{align}
\label{qBI-2}
 \{M_1,M_2\}_q=M_3+m_3,\quad \{M_2,M_3\}_{q}=M_1+m_1,\quad \{M_3,M_1\}_{q}=M_2+m_2,
\end{align}
where the structure constants $m_1$, $m_2$ and $m_3$ read
\begin{align}
\label{Structure-1}
\begin{aligned}
 m_1&=(q^{1/2}+q^{-1/2})\left(\frac{(b-b^{-1})(c-c^{-1})-(a-a^{-1})\Upsilon}{(q-q^{-1})^2}\right),
 \\
 m_2&=(q^{1/2}+q^{-1/2})\left(\frac{(c-c^{-1})(a-a^{-1})-(b-b^{-1})\Upsilon}{(q-q^{-1})^2}\right),
 \\
 m_3&=(q^{1/2}+q^{-1/2})\left(\frac{(a-a^{-1})(b-b^{-1})-(c-c^{-1})\Upsilon}{(q-q^{-1})^2}\right).
 \end{aligned}
\end{align}
The presentation \eqref{qBI-2}, \eqref{Structure-1} is clearly invariant with respect to the simultaneous cyclic permutation of the generators $(M_1,M_2,M_3)$ and arbitrary parameters $(a,b,c)$; it is thus $\mathbb{Z}_3$-symmetric. In this realization, the Casimir operator \eqref{qBi-Cas} of the $q$-deformed Bannai--Ito algebra takes the definite value 
\begin{multline}
\label{Cas-Value}
 \Lambda=
 \left(\frac{(a-a^{-1})(b-b^{-1})(c-c^{-1})\Upsilon}{(q-q^{-1})^2}\right)
 \\
 +\left(\frac{a-a^{-1}}{q-q^{-1}}\right)^2+\left(\frac{b-b^{-1}}{q-q^{-1}}\right)^2+\left(\frac{c-c^{-1}}{q-q^{-1}}\right)^2+\left(\frac{\Upsilon}{q-q^{-1}}\right)^2-\frac{q}{(1+q)^2}.
\end{multline}
In view of the above results, one can conclude that the $q$-deformed Bannai--Ito algebra serves as the covariance algebra for $\mathfrak{osp}_q(1|2)$. 
\begin{remark}
 Let us note that if one takes $a=q^{\alpha}$, $b=q^{\beta}$, $c=q^{\gamma}$ and $\Upsilon=-(q^{\delta}-q^{-\delta})$, the structure constants \eqref{Structure-1} and the Casimir value \eqref{Cas-Value} are identical to the ones arising in the Racah problem for $\mathfrak{osp}_q(1|2)$ \cite{2015_Genest&Vinet&Zhedanov_ArXiv_150105602}.
\end{remark}
\section{Conclusion}
In this Letter, the equitable presentation of the quantum superalgebra $\mathfrak{osp}_q(1|2)$ was displayed. It was observed that $\mathfrak{osp}_q(1|2)$ and $\mathfrak{sl}_q(2)$ are related by the formal transformation $q\rightarrow -q$ and it was established that the $q$-deformed Bannai--Ito algebra arises as the covariance algebra of $\mathfrak{osp}_q(1|2)$.

Under the appropriate reparametrization, the $q$-deformed Bannai--Ito algebra \eqref{qBI-2} with structure constants \eqref{Structure-1} tends to the Bannai--Ito algebra in the $q\rightarrow 1$ limit. Similarly in the limit $q\rightarrow 1$ the quantum superalgebra $\mathfrak{osp}_q(1|2)$ defined in \eqref{OSP} tends to the Lie superalgebra $\mathfrak{osp}(1|2)$ extended by its grade involution, also known as $sl_{-1}(2)$ \cite{2011_Tsujimoto&Vinet&Zhedanov_SIGMA_7_93}. However one observes that the operators of the realization \eqref{XYZ}, \eqref{Realization-2}, \eqref{Realization-1}  do not have a well-defined $q\rightarrow 1$ limit. Consequently, one cannot conclude from the above results that the Bannai--Ito algebra is a covariance algebra of $\mathfrak{osp}(1|2)$. The interesting problem of realizing the Bannai--Ito algebra in terms of the Lie superalgebra $\mathfrak{osp}(1|2)$ thus remains.
\section*{Acknowledgements}
V.X.G. holds an Alexander--Graham--Bell fellowship from the Natural Science and Engineering Research Council of Canada (NSERC). The research of L.V. is supported in part by NSERC. A.Z. wishes to thank the Centre de recherches math\'ematiques for its hospitality.
\footnotesize
\begin{multicols}{2}

\end{multicols}
\end{document}